\documentclass[12pt]{amsart}

\usepackage{hyperref}

\usepackage{amscd,epsfig}

\def\P{{\mathbb P}}
\def\Z{{\mathbb Z}}
\def\Q{{\mathbb Q}}
\def\R{{\mathbb R}}
\def\C{{\mathbb C}}

\def\V{{\mathbb V}}

\def\bH{{\mathbb H}}

\def\A{{\mathcal A}}
\def\F{{\mathcal F}}
\def\H{{\mathcal H}}
\def\O{{\mathcal O}}
\def\U{{\mathcal U}}
\def\cV{{\mathcal V}}

\def\a{{\mathfrak a}}
\def\h{{\mathfrak h}}

\def\D{{\Delta}}

\def\vbar{{\overline{v}}}
\def\Fbar{{\overline{F}}}
\def\Qbar{{\overline{\Q}}}
\def\Btilde{\widetilde{B}}
\def\Ctilde{\widetilde{C}}
\def\Ftilde{\widetilde{\F}}

\def\Vtilde{\widetilde{\cV}}
\def\Xtilde{\widetilde{X}}

\def\Pminus{\P^1 - \{0,1,\infty\}}
\def\PminusZ{\P_\Z^1 - \{0,1,\infty\}}
\def\Pminusp{(\Pminus)}
\def\v{{\vec{v}}}
\def\b{{\vec{01}}}
\def\c{{\vec{10}}}

\def\Ql{{\Q_\ell}}

\def\SL{{\mathrm{SL}}}
\def\Sp{{\mathrm{Sp}}}
\def\Gm{{\mathbb{G}_m}}

\def\dot{\bullet}
\def\blank{{\underline{\phantom{x}}}}
\def\comp{~\widehat{\!}{\;}}

\def\MZN{\mathrm{MZN}}
\def\DR{\mathrm{DR}}

\def\lb{\langle\langle}
\def\rb{\rangle\rangle}

\def\limproj#1{\lim_{\stackrel{\longleftarrow}{#1}}}

\newcommand\id{\operatorname{id}}
\newcommand\Aut{\operatorname{Aut}}
\newcommand\End{\operatorname{End}}
\newcommand\Der{\operatorname{Der}}
\newcommand\Hom{\operatorname{Hom}}
\newcommand\Ext{\operatorname{Ext}}
\newcommand\Jac{\operatorname{Jac}}
\newcommand\Spec{\operatorname{Spec}}
\newcommand\im{\operatorname{im}}
\newcommand\Gr{\operatorname{Gr}}
\renewcommand\Im{\operatorname{Im}}
\newcommand\Gal{\operatorname{Gal}}

\newtheorem{theorem}{Theorem}		%[section]
\newtheorem{lemma}[theorem]{Lemma}

\newtheorem{corollary}[theorem]{Corollary}

\theoremstyle{definition}
\newtheorem{definition}[theorem]{Definition}
\newtheorem{example}[theorem]{Example}

\theoremstyle{remark}
\newtheorem{remark}[theorem]{Remark}

\begin{document}

\title{Periods of Limit Mixed Hodge Structures}
\author{Richard Hain}

\dedicatory{To Wilfred Schmid on the occasion of this 60th birthday}

\maketitle

\section{Introduction}

The first goal of this paper is to explain some important results of Wilfred
Schmid from his fundamental paper \cite{schmid} in which he proves very general
results which govern the behaviour of the {\it periods} of a of smooth
projective variety $X_t$ as it degenerates to a singular variety. As has been
known since classical times, the periods of a smooth projective variety
sometimes contain significant information about the geometry of the variety,
such as in the case of curves where the periods determine the curve. Likewise,
information about the asymptotic behaviour of the periods of a variety as it
degenerates sometimes contain significant information about the degeneration
and the singular fiber. For example, the {\it Hodge norm} estimates, which are
established in \cite{schmid} and \cite{cat-kap-schmid:sl2}, describe the
asymptotics of the Hodge norm of a cohomology class as the variety degenerates
in terms of its monodromy. They are an essential ingredient in the study of the
$L_2$ cohomology of smooth varieties with coefficients in a variation of Hodge
structure \cite{zucker,cat-kap-schmid:ell2}.

A second goal is to give some idea of how geometric and arithmetic information
can be extracted from the limit periods, both in the geometric case and in the
case of the limits of the mixed Hodge structures on fundamental groups of
curves.

To get oriented, recall that if $X$ is a compact Riemann surface of genus $g$,
then for each choice  of a symplectic (w.r.t.\ the intersection form) basis
$a_1,\dots,a_g,b_1,\dots,b_g$ of $H_1(X,\Z)$, there is a basis $w_1,\dots,w_g$
of the holomorphic differentials $H^0(X,\Omega^1)$ such that
$$
\int_{a_j} w_k = \delta_{jk}.
$$
The $g\times g$ matrix
$$
\Omega := \bigg(\int_{b_j} w_k\bigg)
$$
is called the {\it period matrix} of $X$. The classical Riemann bilinear
relations assert that $\Omega$ is symmetric and that it has positive definite
imaginary part. By the classical Torelli Theorem \cite{griffiths-harris}, the
period matrix determines the Riemann surface up to isomorphism.

A more conceptual way to view the period matrix is to note that the augmented
period matrix $(I|\Omega)$ is the set of coordinates of the $g$-plane
$H^0(X,\Omega^1)$ in the grassmannian of $g$ planes in the $2g$-dimensional
vector space $H^1(X,\C)$ with respect to the integral basis $a_1,\dots,b_g$ of
$H_1(X,\Z)$. The $g$-plane $H^0(C,\Omega^1)$ lies in the closed submanifold of
the grassmannian consisting of $g$-planes isotropic with respect to the cup
product.

\section{Quick Review of Hodge Theory}

\subsection{The Hodge Theorem}

For any complex manifold, we can define
$$
H^{p,q}(X) :=
\frac{\text{closed $(p,q)$-forms}}{\text{$(p,q)$-forms that are exact}}.
$$
This is simply the subspace of $H^{p+q}(X,\C)$ consisting of classes that can
be represented by a closed $(p,q)$-form.

For a general complex manifold, these do not give a decomposition of the
cohomology of $X$. However, if $X$ is compact K\"ahler (for example, a smooth
projective variety), then the natural mapping
$$
\bigoplus_{p+q = k} H^{p,q}(X) \to H^k(X,\C)
$$
is an isomorphism. Denote the class in $H^{1,1}(X)$ of the K\"ahler form by
$w$. The Hard Lefschetz Theorem states that if $\dim X = n$ and $k\ge 0$, then
$$
\blank \wedge w^k : H^{n-k}(X,\C) \to H^{n+k}(X,\C)
$$
is an isomorphism. One can then define the {\it primitive} cohomology in
degree $k$ by
$$
PH^k(X,\C) =
\ker\big\{\blank\wedge w^{k+1} : H^{n-k}(X,\C) \to H^{n+k+2}(X,\C)\big\}.
$$
The Hodge decomposition of $H^k(X)$ restrictions to give one of the primitive
cohomology:
$$
PH^k(X,\C) \cong \bigoplus_{p+q=k} PH^{p,q}(X)
$$
where $PH^{p,q}(X) = H^{p,q}(X) \cap PH^{p+q}(X,\C)$.

Define a form $Q : H^{n-k}(X) \otimes H^{n-k}(X) \to \C$ by 
$$
Q(\xi,\eta) = \int_X \xi \wedge \eta \wedge w^k.
$$
This form is always defined over $\R$ as $w$ is a real cohomology class, and
over $\Z$ when $X$ is smooth projective and $w$ is the class of a hyperplane
section.

The Riemann-Hodge bilinear relations generalize the classical Riemann bilinear
relations for a compact Riemann surface. They state that
\begin{enumerate}
\item $Q$ vanishes on $PH^{p,q}(X)\otimes PH^{r,s}(X)$ unless $p=s$ and $q=r$;
\item if $\xi \in PH^{p,q}(X)$ is non-zero, then
$$
i^{p-q}(-1)^{k(k-1)/2} Q(\xi,\overline{\xi}) > 0,
$$
where $k=p+q$.
\end{enumerate}
In degree 2, the Riemann-Hodge bilinear relations imply the Hodge Index
Theorem.

The primitive cohomology $P^kH^{p,q}(X)$ together with the form
$S:=(-1)^{k(k-1)/2}Q$ is the prototypical example of a {\it polarized Hodge
structure of weight} $k$.

\begin{definition}
A {\it Hodge structure $V$ of weight $k$} consists of
\begin{enumerate}
\item a finitely generated abelian group $V_\Z$ and
\item a bigrading
$$
V_\C = \bigoplus_{p+q=k} V^{p,q}
$$
of its complexification, which satisfies $\overline{V^{p,q}} = V^{q,p}$.
\end{enumerate}
A {\it polarized Hodge structure of weight $k\in\Z$} is a Hodge structure $V$
of weight $k$ together with a $(-1)^k$ symmetric bilinear form
$$
S : V_\Z \otimes V_\Z \to \Z
$$
that satisfies:
\begin{enumerate}
\item $S(V^{p,q},V^{r,s}) = 0$ unless $p=s$ and $q=r$,
\item if $v \in V^{p,q}$ is non-zero, then $i^{p-q}S(v,\vbar) > 0$.
\end{enumerate}
\end{definition}

The Hodge norm of $v \in V$ is defined by
$$
\|v\|^2 = S(Cv,\vbar)
$$
where $C$ is the linear operator whose restriction to $V^{p,q}$ is $i^{p-q}$.

\section{Periods}

The {\it Hodge filtration} of a Hodge structure $V$ of weight $k$ is the
decreasing filtration
$$
 \cdots \supseteq F^pV \supseteq F^{p+1}V \supseteq F^{p+2} V \supseteq \cdots
$$
defined by
$$
F^p V := \bigoplus_{s\ge p} V^{s,k-s}.
$$
One can recover the bigrading from the Hodge filtration and the real structure
by
$$
V^{p,k-p} = F^p V \cap \Fbar^{k-p},
$$
where $\Fbar^\dot$ is the complex conjugate of the Hodge filtration.

The reason for working with the Hodge filtration rather than the bigrading is
the observation of Griffiths \cite{griffiths} that the Hodge filtration varies
holomorphically\footnote{With respect to the locally constant structure coming
from the lattice of integral cohomology.} in families, whereas the $(p,q)$
pieces generally do not. Indeed, if $H^{p,q}$ varies holomorphically, then
$H^{q,p}$ varies anti-holomorphically. So if both $H^{p,q}$ and $H^{q,p}$ vary
holomorphically, then both are locally constant.

To be more precise, suppose that $f:X \to B$ is a family of compact K\"ahler
manifolds. This means that $f$ is a proper holomorphic mapping, that $X$ is
K\"ahler and that each fiber $X_b$ is smooth and has the induced K\"ahler
structure. One has the local system (i.e., locally constant sheaf)
$$
\bH^k := R^k f_\ast \Z = \{H^k(X_b,\Z)\}_{b\in B}
$$
on $B$. We can complexify this to obtain a flat holomorphic vector bundle
$$
\H^k := \bH^k \otimes_\Z \O_B
$$
over $B$. Denote the  natural flat connection on $\H^k$ by $\nabla$. Let
$$
\F^p := \big\{F^pH^k(X_b)\big\}_{b\in B} \subseteq \H^k.
$$
Fundamental results of Griffiths \cite{griffiths} assert that the $\F^p$ are
holomorphic sub-bundles of $\H^k$ and that if one differentiates a section of
$\F^p$ along a holomorphic vector field, then the result lies in $\F^{p-1}$:
\begin{equation}
\label{inf_per_reln}
\nabla : \F^p \to \F^{p-1}\otimes \Omega^1_B.
\end{equation}
This is the prototypical example of a {\it variation of Hodge structure}.
Replacing $\bH^k$ by the bundle of primitive cohomology
$$
P\bH^k := \{PH^k(X_b) \}_{b \in B}
$$
yields the prototypical example of a {\it polarized variation of Hodge
structure} --- a variation of Hodge structure with an inner product $S$ which
is parallel with respect to the connection and which polarizes each fiber.

\begin{definition}
A {\it variation of Hodge structure of weight $k$} over a complex manifold $B$
is a $\Z$-local system $\V$ together with a flag
$$
\cdots \supseteq \F^p \supseteq \F^{p+1} \supseteq \cdots
$$
of holomorphic sub-bundles of the flat bundle $\cV := \V\otimes_\Z \O_B$ which
satisfy:
\begin{enumerate}
\item
$\nabla : \F^p \to \F^{p-1}  \otimes \Omega^1_B$,
\item Each fiber, endowed with the induced Hodge filtration, is a Hodge
structure of weight $k$.
\end{enumerate}
The variation is {\it polarized} if it has a $(-1)^k$-symmetric inner product
which is parallel with respect to the flat structure and polarizes the Hodge
structure on each fiber.
\end{definition}

Given a polarized variation of Hodge structure $\V$ of weight $k$ over a base
$B$, one has a period  mapping
$$
\Btilde \to
\left\{\text{set of flags $\{F^p\}$ that satisfy $S(F^p,F^{k-p+1})=0$} \right\}
$$
defined on the universal covering of the base. The flags should be considered
to live in some fixed fiber of $\V \to B$; other fibers being compared with
this one using the flat structure. The flags should have the ``same shape'' as
the Hodge filtration on this reference fiber.

Griffiths' results assert that this mapping is holomorphic and lies in the open
subset $D$ of flags that satisfy
$$
i^{p-q}S(v,\vbar) > 0
$$
whenever $v$ is a non-zero element of $F^p\cap \Fbar^{k-p}$ and $p+q=k$.

When the variation is $P\bH^k$ this mapping really is a period mapping in the
sense that the coordinates in the flag manifold of the image of $b\in \Btilde$
will be given by periods over integral cycles of a basis of $PH^k(X_b)$ adapted
to the Hodge filtration.

At this point it is useful to revisit the period mapping for curves.

\begin{example}
Fix a genus $g \ge 1$. To know the Hodge filtration on a compact Riemann
surface of genus $g$, one simply has to know the location of the holomorphic
differentials $F:=H^{1,0}(X_b)$ in $H^1(X,\C)$ with respect to an integral
symplectic basis. The polarization $S$ on $H^1(X_b)$ is the one given by cup
product. The first Riemann bilinear relation says that $S(F,F) = 0$.

Choose a base point $b_o \in B$. Let $H = H^1(X_{b_o})$. The image of the
period mapping is contained in the closed submanifold $Y$ of the grassmannian
of $g$ planes in $H_\C$ consisting of those $g$-planes $F$ that satisfy $S(F,F)
= 0$. It is a complex manifold of dimension $g(g+1)/2$. The period mapping
$\Btilde \to Y$ is holomorphic, and has image contained in the open subset
$$
U = \{F \in Y : i S(v,\vbar) > 0 \text{ if } v\neq 0, v \in F\}
$$
of $Y$; it is a complex manifold. One can check that
\begin{align*}
U
&= \text{ the symmetric space of $\Sp_g(\R)$} \cr
&= \Sp_g(\R)/U(g) \cr
&= \text{ Siegel's upper half plane $\h_g$ of rank $g$} \cr
&= \{g\times g \text{ symmetric complex matrices $\Omega$ with
$\Im \Omega >0$}\}.
\end{align*}
The period map takes $X_b$ to its period matrix $\Omega(b)$, described in the
Introduction. That it is holomorphic is equivalent to the statement that
$\Omega(b)$ depends holomorphically on $b\in B$.

The monodromy representation of the variation is a homomorphism
$$
\pi_1(B,b_o) \to \Sp_g(\Z).
$$
The period mapping descends to the mapping
$$
B \to \Sp_g(\Z) \backslash \h_g = \A_g
$$
to the moduli space of principally polarized abelian varieties that takes $b
\in B$ to the jacobian of $X_b$.
\end{example}

\section{The Limit Hodge Filtration}
\label{lim_mhs}

Schmid's work concerns the behaviour of the period mapping of a polarized
variation of Hodge structure as the the Hodge structures degenerate. In the
geometric case, this happens when the varieties $X_b$ degenerate to a singular
variety. To keep things simple I'll only discuss the local case in one
variable, which he considers in \cite{schmid}. The several variable case is
more subtle and is worked out in \cite{schmid} and in \cite{cat-kap-schmid:sl2}.

Suppose that $\V \to \D^\ast$ is a polarized variation of Hodge structure
of weight $k$ over the punctured disk. The geometrically inclined may prefer
to think of this as the local system of primitive cohomology of degree $k$
associated to a proper family $X \to \D$ over the disk, where $X$ is smooth
and K\"ahler and where each fiber $X_t$ is smooth when $t\neq 0$.

Let
$$
T : V_{t_o} \to V_{t_o}
$$
be the monodromy operator associated to the positive generator of
$\pi_1(\D^\ast,t_o)$. Landman in the geometric case, and Borel in general,
(reproved by Schmid \cite[p.~245]{schmid}) proved that the eigenvalues of $T$
are roots of unity, so that $T$ is quasi unipotent; that is, there exist
integers $e$ and $m$ such that
$$
(T^e-I)^m = 0.
$$
By pulling $\V$ back along the $e$-fold covering $t\mapsto t^e$ of $\D^\ast \to
\D^\ast$ we may assume that $T$ is unipotent. The logarithm of $T$ is then
defined by the usual power series. Set
$$
N = \log T/2\pi i.
$$
This is nilpotent and preserves the polarization in the sense that
\begin{equation}
\label{inf_sym}
S(Nx,y) + S(x,Ny) = 0.
\end{equation}
We shall view it as a flat section of the bundle $\End\V$, so that it acts on
all fibers, not just the one over $t_o$.

Denote the flat holomorphic vector bundle $\V\otimes_\Z \O_{\D^\ast}$ by $\cV$
and its connection by $\nabla$.

Note that, unless $T$ has finite order, neither $T$ nor $N$ is a morphism of
Hodge structures $V_{t_o} \to V_{t_o}$. One of Schmid's results, which is
explained below, is that $N$ is compatible with the  limit Hodge filtration and
acts as an endomorphism of the ``limit mixed Hodge structure'', which allows
one to relate Hodge theory and monodromy.

\subsection{Deligne's Canonical Extension}
In order to discuss the limit of the Hodge structure $V_t$ as $t\to 0$, we need
to extend the flat bundle $\cV \to \D^\ast$ to a holomorphic vector bundle
$\Vtilde \to \D$ over the entire disk. The way to do this is to use Deligne's
{\it canonical extension} \cite{deligne:ode}.

First note that any trivialization of $\cV$ over $\D^\ast$ gives an extension
of $\cV$ to $\D$; the holomorphic sections of the extended bundle are those
whose restriction to $\D^\ast$ are of the form
$$
\sum_j f_j \phi_j
$$
where each $f_j$ is holomorphic on $\D$, and $\phi_1,\dots , \phi_m$ is a
holomorphic frame that gives the trivialization over $\D^\ast$. Two framings of
$\cV$ over $\D^\ast$ give the same extension if and only if they differ by a
holomorphic mapping $g : \D \to GL_m(\C)$. It follows that the set of
extensions of $\cV$ to $\D$ can be identified with the homogeneous space
$$
GL_m(\O(\D^\ast))/GL_m(\O(\D)),
$$
which is large.

To construct Deligne's extension, we choose multivalued flat sections
$\phi_1(t),\dots,\phi_m(t)$ defined on $\D^\ast$. The extension is constructed
by regularizing these. Specifically, set
$$
e_j(t) = \phi_j(t) t^{-N},
$$
where for a square matrix $A$, $t^A$ is defined to be $\exp(A \log t)$. The
$e_j(t)$ comprise a single valued framing of $\cV$ over $\D^\ast$ as each
$\phi_j(t)$ changes to
$$
 \phi_j(t) \cdot T
$$
when analytically continued around the circle, whereas $t^{-N}$ changes to
$$
\exp\big((-\log t - 2\pi i)N\big) = e^{-2 \pi i N} t^{-N}
= T^{-1} t^{-N}
$$
when continued around the origin. The framing $e_j$ determines Deligne's
extension $\Vtilde$ of $\cV$ to all of $\D$. With respect to Deligne's
framing $e_j(t)$, the connection is
$$
\nabla = d - N \frac{dt}{t}.
$$
In other words, the connection on $\Vtilde$ is meromorphic with a simple pole
at the origin and the residue of the connection at $t=0$ is $-N$. In
particular, the extended connection has a regular singular point at the origin
with nilpotent residue. This property characterizes Deligne's extension up to
isomorphism. This fact, proved by Deligne in more generality in
\cite{deligne:ode}, is a consequence of the classical result about differential
equations stated below, a proof of which can be found in
\cite[Chapt.~II]{wasow}.

\begin{lemma}
Suppose that $A : \D \to \End V$ is holomorphic. If no two eigenvalues of
$A_0 := A(0)$ differ by a non-zero integer (e.g., if $A_0$ is nilpotent), then
there is a unique holomorphic function $P : \D \to \Aut V$ with $P(0) = \id_V$
such that each (local) solution $v : \D^\ast \to V$ of the differential
equation
$$
t v'(t) = v(t)A(t)
$$
is of the form
$$
v(t) = v_0\, t^{A_0} P(t).
$$
\end{lemma}

Note that the monodromy operator $T$ is:
$$
v(t) = v_0\,t^{A_0}P(t) \mapsto v_0\, t^{A_0} e^{2\pi i A_0} P(t)
=  v(t) \exp\big({2\pi i P(t)^{-1}A_0\, P(t)}\big).
$$
The operator $\log T/2\pi i$ on the fiber over $t$ is thus right multiplication
by
$$
N(t) = P(t)^{-1} A_0\, P(t)\ \in \Aut V.
$$
The Deligne trivializing sections are therefore of the form
$$
v(t) t^{-N(t)} = v(t) P(t)^{-1} t^{-A_0} P(t) =  v_0\, P(t),
$$
which converges to $v_0$ as $t \to 0$. This also implies that the function
$P(t)$ is the difference between the Deligne framing and the given
trivialization.

\begin{corollary}
\label{limit}
If $A_0$ is nilpotent, then
$$
\lim_{t\to 0} v(t) t^{-A_0} = \lim_{t\to 0} v(t) t^{-N(t)} = v_0
$$
where the limit is taken along any angular ray.
\end{corollary}

\begin{proof}
The second equality was proved above. As for the first, we have
$$
v(t) t^{-A_0} = v_0\, t^{A_0} P(t) t^{-A_0}.
$$
If $A_0^{k+1} = 0$, then there is a constant $C$ such that
$$
\|t^{A_0}\| \text{ and } \|t^{-A_0}\| \le C \big(\log1/|t|\big)^k
$$
when $0<|t|\le R$, for some $R>0$. Writing $P(t) = I + \sum_{n\ge 1} P_n t^n$,
we have
$$
\|t^{A_0} P(t) t^{-A_0} - I\| \le
2C|t|\big(\log1/|t|\big)^k \sum_{n=1}^\infty \|P_n\| |t|^{n-1}
$$
which goes to zero along each angular ray.
\end{proof}

This formula will prove useful in Section~\ref{mhs-p1} when computing the period
of the path $[0,1]$ in the limit MHS on the unipotent completion in the set of
paths in $\P^1 - \{0,1,\infty\}$ from $0$ to $1$.

\begin{remark}
Note that $v(t)t^{-A_0}$ is not, in general, single-valued on the punctured
disk even though its limit as $t\to 0$ along each radial ray exists.
\end{remark}

\subsection{The Nilpotent Orbit Theorem}

Schmid's first main result is the Nilpotent Orbit Theorem. We shall state it as
two separate results.

\begin{theorem}[Schmid]
The Hodge sub-bundles $\F^p$ of $\cV$ extend to holomorphic sub-bundles
$\Ftilde^p$of Deligne's canonical extension $\Vtilde \to \D$. That is, both
$\Ftilde^p$ and $\Vtilde/\Ftilde^p$ are vector bundles.
\end{theorem}

This shows that, when taken appropriately, the limit of the Hodge filtration
exists. The limit Hodge filtration $\Ftilde^\dot_0$ is a filtration of the
fiber $V_0$ of the canonical extension $\Vtilde$ over 0.

In the geometric case, the extended Hodge bundles can be constructed
explicitly. This was done by Steenbrink in \cite{steenbrink}.

So far we have constructed a complex vector space $V_0$ and a Hodge filtration
$F^\dot V_0 := \Ftilde^p_0$ on it. The polarization also extends to the central
fiber as at $t^N$ preserves the inner product since $N$ preserves it
infinitesimally by (\ref{inf_sym}).

We can also construct an integral form $V_{\Z}$ of $V_0$. If the $\phi_j$ are
integral, the integral lattice $V_\Z$ is simply the $\Z$-linear span of the
$e_j(0)$. It is important to note, however, that except when the monodromy is
trivial, $(V_\Z,F^\dot_0)$ is not a Hodge structure even though it is a limit
of Hodge structures.

In general, the integral structure depends on the holomorphic parameter $t$
chosen for the disk. This dependence is quite transparent, and is a
straightforward consequence of Corollary~\ref{limit}.

\begin{lemma}
The lattice $V_\Z$ in $V_0$ depends only on the parameter $t$ to first order,
and hence only on the tangent vector $\partial/\partial t$. The lattice
corresponding to the tangent vector $\lambda \partial/\partial t$, where
$\lambda$ is a sufficiently small non-zero complex number, is
$\lambda^N \cdot V_\Z$. \qed
\end{lemma}

This extends to a ``nilpotent orbit'' over $\C^\ast$. Its fiber over $t \in
\C^\ast$ consists of:
\begin{enumerate}
\item the complex vector space $V_0$,
\item the limit Hodge filtration $F^\dot_0$ on $V_0$,
\item the lattice $t^N V_\Z$ in $V_0$.
\end{enumerate}
It  has a flat connection given by the limit connection
$$
\nabla_0 = d - N\frac{dt}{t}.
$$
with respect to the constant trivialization given by $V_0$. It satisfies the
Griffiths infinitesimal period relation (\ref{inf_per_reln}).

The Deligne frame gives a connection preserving isomorphism
$$
\begin{CD}
V_0 \times \D @>>> \Vtilde \cr
@VVV @VVV \cr
\D @= \D
\end{CD}
$$
between $\Vtilde$ and the restriction of the nilpotent orbit to a neighbourhood
of the origin which is the identity on the fiber over the origin. The second
part of Schmid's Nilpotent Orbit Theorem can be stated as follows.

\begin{theorem}[Schmid]
The nilpotent orbit satisfies Griffiths infinitesimal period relation
(\ref{inf_per_reln}). There is an $\epsilon > 0$ such that if $0 < |t| <
\epsilon$, then the fiber of the nilpotent orbit over $t$ is a Hodge structure
with the same Hodge numbers as those of the original  variation. Moreover,
there is a holomorphic mapping $g : \D \to \Aut(V_0,S)$ such that $g(0)$ is the
identity and $g$ carries the Hodge bundle of the nilpotent orbit onto the Hodge
bundle of $\Vtilde$.
\end{theorem}

So this gives us our first precise information about the asymptotics of periods
of polarized variations of Hodge structure. Basically, the crudest model of a
degeneration of Hodge structure is a nilpotent orbit. General degenerations are
obtained from these by perturbing the Hodge filtration.

It is useful to think of the nilpotent orbit as being defined on the punctured
tangent space $T_0\D -\{0\}$ of $\D$ at $t=0$. Then in a sufficiently small
neighbourhood of the origin, the nilpotent orbit is a variation of Hodge
structure that approximates the original variation. The nilpotent orbit can be
thought of as the restriction of the canonical extension of the variation $\V$
to $\Spec \C[t]/(t^2)$.

\section{The $\SL_2$-orbit Theorem}

Schmid's second main theorem, and the deeper of the two, is the $\SL_2$-orbit
Theorem. It gives more precise information about how nilpotent orbits
approximate degenerations of polarized variations of Hodge structure and is
analogous to the Hard Lefschetz Theorem for the cohomology of a smooth
projective variety. In fact, the $\SL_2$-orbit Theorem corresponds to the Hard
Lefschetz Theorem for Calabi-Yau manifolds under mirror symmetry.\footnote{As
Eduardo Cattani pointed out to me (cf.\ \cite{cattani}), a precursor of this
already appears in the papers \cite{cat-kap-schmid:sl2,cat-kap-schmid:ell2} of
Cattani-Kaplan-Schmid. They prove a general result, a special case of which
constructs out of $(H^\dot(X),F^\dot)$, the cohomology and Hodge filtration of
a compact K\"ahler manifold $X$, a polarized variation of Hodge structure over
the complexified K\"ahler cone of $X$. This can be restricted to the the
$(p,p)$ part of $H^\dot(X)$, and this is what should correspond to the
polarized variation of Hodge structure of the degeneration of the mirror, as
was pointed out to David Morrison by Cattani in 1994. Deligne
\cite{deligne:mirror} has also studied certain nilpotent orbits arising in
mirror symmetry.} It also leads naturally to mixed Hodge structures, which,
like Hodge structures, have a lattice and Hodge filtration, but also have a
{\it weight filtration}. In the case of a limit of Hodge structures, the weight
filtration is constructed from the monodromy logarithm.

\subsection{The Monodromy Weight Filtration}
\label{monod_wt_filt}

Suppose that $V$ is a finite dimensional vector space over a field of
characteristic zero, and suppose that $N$ is a nilpotent endomorphism of $V$.
Then there is a unique filtration $W(N)_\dot$
$$
\cdots \subseteq W_{n}(N) \subseteq W_{n+1}(N) \subseteq \cdots
$$
of $V$ such that
\begin{enumerate}
\item $N(W_n(N)) \subseteq W_{n-2}(N)$ and
\item $N^k : W_k/W_{k-1} \to W_{-k}/W_{-k-1}$ is an isomorphism.
\end{enumerate}

This is an easy exercise using the Jordan canonical form of $N$. It suffices
to prove the result for a single Jordan block, where the result is evident.

In our case, the monodromy logarithm acts on $V_0$. We therefore have the
monodromy weight filtration $W_\dot(N)$ of $V_0$. In the present situation,
we shift the filtration by the weight $k$ of the variation --- define:
$$
W_n V_0 = W_{k-n}(N).
$$

Our central fiber now has the Hodge filtration $F^\dot V_0$, the nilpotent
orbit $t^N V_\Z$ of integral structures, and the monodromy weight filtration
$W_\dot$. One of the main consequences of the $\SL_2$-Orbit Theorem is:

\begin{theorem}[Schmid]
For each $t \in \C^\ast$, the collection
$$
(t^N V_\Z,F^\dot_0,W_\dot)
$$
is a mixed Hodge structure. Moreover, $N : V_0 \to V_0$ is defined
over $\Q$ and is is a morphism of $\Q$-mixed Hodge structure of type $(-1,-1)$.
\end{theorem}

So, even if one is interested only in compact K\"ahler manifolds, one is 
inexorably lead to mixed Hodge structures when one studies their degeneration.

\begin{definition}[Deligne]
A {\it mixed Hodge structure} $V$ consists of a finitely generated abelian
group $V_\Z$, an increasing filtration
$$
\cdots \subseteq W_{m-1}V_\Q \subseteq W_m V_\Q \subseteq W_{m+1}V_\Q \subseteq
\cdots
$$
(called the {\it weight filtration}) of $V_\Q := V_\Z \otimes_\Z \Q$, and a
decreasing filtration
$$
\cdots \supseteq F^{p-1} V_\C \supseteq F^p V_\C \supseteq F^{p+1} V_\C
\supseteq \cdots
$$
(called the {\it Hodge filtration}) of $V_\C := V_\Z\otimes_\Z \C$. These are
required to satisfy the condition that, for each $m$,
$$
\Gr^W_m V := W_{m}V/W_{m-1}V
$$
with the induced Hodge filtration
$$
F^p \Gr^W_m V := \im\{ F^p V \cap W_m V \to \Gr^W_m V\}
$$
is a Hodge structure of weight $m$.
\end{definition}

Deligne \cite{deligne:II,deligne:III} proved that the cohomology groups of
every complex algebraic variety (not necessarily smooth or compact) have a
mixed Hodge structure (MHS) that is natural with respect to morphisms.

This definition takes a while to digest, and it may be helpful to note that, up
to some subtle issues of torsion in $\V_\Z$, mixed Hodge structures can be
constructed as follows.

Let $(H^m)_{m\in \Z}$ be a collection of Hodge structures where $H^m$ has
weight $m$ and all but a finite number of the $H^m$ are zero. Then their direct
sum
$$
H  = \bigoplus_m H^m
$$
is a mixed Hodge structure with the direct sum Hodge and weight filtrations:
$$
F^p H := \bigoplus_k F^p H^m,\quad W_m H = \bigoplus_{k \le m} H^k.
$$
Such a MHS is said to be {\it split} as it is the direct sum of Hodge
structures. Very few of the MHSs encountered in life are split. Up to the
torsion issues alluded to above, all MHSs $V$ with $\Gr^W_m V \cong H^m$ can be
constructed as follows.

Let
$$
G = \{\phi \in \Aut H_\C \text{ that preserve $W_\dot$ and act trivially
on all $\Gr^W_m H$}\}.
$$
This is a unipotent group. For each element of $g$ of $G$, we can define $V(g)$
to be the MHS with $V_\C = H_\C$, the same Hodge and weight filtrations (over
$\C$), and lattice $g\cdot H_\Z \hookrightarrow H_\C$. The MHS obtained in this
way, depends only on the double coset of $g^{-1}$ in
$$
G_\Z \backslash G / F^0 G
$$
where
$$
G_\Z = G\cap \Aut V_\Z \text{ and } F^0 G = \{ \phi \in G : \phi\text{
preserves } F^\dot\}.
$$
In fact, this double coset space is the moduli space of MHSs $V$ together
with an isomorphism $\Gr^W_m V \cong H^m$.

\begin{example}
Let $A$ and $B$ be Hodge structures of weights $a$ and $b$, respectively, where
$b<a$. The set of all mixed Hodge structures with weight graded quotients
isomorphic to $A$ and $B$, together with these {\it framings} is
$$
\Hom(A_\C,B_\C)/\big(\Hom(A_\Z,B_\Z) + F^0 \Hom(A,B)\big)
$$
as the group $G$ in this case is $\Hom(A_\C,B_\C)$. This gives the computation
of $\Ext_\mathrm{Hodge}^1(A,B)$ in the category of MHSs given by Carlson in
\cite{carlson}.
\end{example}

\subsection{The Hodge Norm Estimates}

Another important consequence of the $\SL_2$-orbit Theorem is the Hodge norm
estimates, which are important in applications to $L_2$-cohomology of algebraic
varieties with coefficients in a polarized variation of Hodge structure ---
cf.\ \cite{zucker} in the one-variable case, and \cite{cat-kap-schmid:ell2} and
\cite{looijenga} in general.

First note that the monodromy weight filtration is defined on the cohomology of
each fiber and is preserved by the connection.

\begin{theorem}[Schmid]
If $\V$ is a polarized variation of Hodge structure of weight $k$ over the
punctured disk $\D^\ast$, then a flat section $v(t)$ lies in $W_m$ if and
only if the square of its Hodge norm satisfies
$$
S(Cv(t),\vbar(t)) = O\bigr(\big(\log(1/|t|)\big)^{m-k}\bigr)
$$
when $t$ goes to 0 along a radial ray. Here $C : V_t \to V_t$ is the operator
that is multiplication by $i^{p-q}$ on $V_t^{p,q}$.
\end{theorem}

\subsection{Topology and Geometry of the Limit MHS}
Suppose that $f:X \to \D$ is a proper holomorphic mapping from a K\"ahler
manifold to the disk. Suppose that the fiber $X_t$ over $t\in \D$ is smooth
whenever $t\neq 0$ and that $X_0$ is a reduced divisor with normal crossings.
This implies that the monodromy operator $T \in \Aut H^m(X_t)$ is unipotent.

The local system
$$
\V := R^m \big(f|_{\D^\ast}\big)_\ast\Z
$$
underlies a polarized variation of Hodge structure of weight $m$ over
$\D^\ast$. By Schmid's Theorems, for each choice of a non-zero tangent vector
$\v$ of $\D$ at $0$, there is a limit mixed Hodge structure of $\V$, which one
can think of as a MHS on the cohomology group $H^m(X_\v)$ of the first order
deformation $X_\v$ of $X_0$ given by $\v$.

Basic results about the compatibility of this MHS with others associated to the
degeneration follow from the geometric constructions of the limit MHS given by
Clemens \cite{clemens} and Steenbrink \cite{steenbrink}. These compatibilities
aid in extracting geometric and topological information from $H^\dot(X_\v)$.

To begin to understand the limit MHS, $H^m(X_\v)$, it is helpful to think of
$X_v$ as being built out of $X_0$. In order to explain this, we denote the
normalization of $X_0$ by $\nu : \Xtilde_0 \to X_0$. Denote the inverse image
under $\nu : \Xtilde_0 \to X_0$ of the singular locus of $X_0$ by $D$. It is a
normal crossings divisor. Set
$$
X_0' = \text{ real oriented blowup of $\Xtilde_0$ along  $D$}
$$
which is a manifold with corners. One can think of each $X_\v$ as being
obtained from $X_0'$ by gluing, the data for which is given by the
combinatorics of $D$, its local defining equation of $X_0$, and by $\v/|\v|$.
There are therefore mappings
$$
\Xtilde_0 - D \hookrightarrow X_0' \to X_\v \to X_0
$$
which induce homomorphisms
$$
\begin{CD}
H^m(X_0) @>\beta>> H^m(X_v) @>\alpha>> H^m(\Xtilde_0 - D).
\end{CD}
$$

\begin{theorem}
The mappings $\alpha$ and $\beta$ are both morphisms of MHS.
\end{theorem}

The statement that $\beta$ is a morphism is proved by both Clemens and
Steenbrink. It is closely related to the {\it Local Invariant Cycle Theorem},
which is stated below. That $\alpha$ is a morphism follows from Steenbrink's
construction of $H^\dot(X_\v)$ in \cite{steenbrink}.

For each $t\in \D$, there is a restriction mapping
$$
j_t^\ast : H^m(X,\Q) \to H^m(X_t,\Q)
$$
which corresponds, via duality, to intersection with the fiber
$$
H_n(X,X|_{\partial\D},\Q) \to H_{n-2}(X_t,\Q)
$$
and to the mapping $\beta$ above after taking into account the fact that the
inclusion $X_0 \hookrightarrow X$ is a deformation retraction.

The following result was first established in the $\ell$-adic case by Deligne
and then in the complex analytic case by Clemens \cite{clemens} and Steenbrink
\cite{steenbrink}. Although not expressed in terms of limit Hodge theory, it
is a non-trivial application of its existence and construction.

\begin{theorem}[Local Invariant Cycle Theorem]
If $t\neq 0$, then the image of $j_t^\ast$ is the space
$$
\ker\big\{(T-I) : H^m(X_t,\Q) \to H^m(X_t,\Q)\big\}
$$
of invariant cohomology classes.
\end{theorem}

\subsection{Degenerations of Curves}

In this section, we give a relatively simple example to illustrate how
geometric information can be extracted from the limit mixed Hodge structure. To
get the full power from this theory, one needs to combine Schmid's results with
the complementary results of Clemens \cite{clemens} and the explicit
constructions of Steenbrink \cite{steenbrink} surveyed in the previous
paragraph.

Suppose that $C \to \D$ is a stable degeneration of compact Riemann
surfaces of genus $g$. The total space $C$ is assumed to be smooth and the
fiber $C_t$ over $t$ is assumed to be smooth when $t\neq 0$. The central
fiber is assumed to be reduced and stable (i.e., its automorphism group
is finite.)

Let $B$ be the set of the homology classes of the vanishing cycles. Then,
using the Picard-Lefschetz formula, one sees that the monodromy is unipotent
and satisfies $(T-I)^2 = 0$. This implies that $N = T-I$ and that the monodromy
weight filtration has length 3:
$$
0 \subseteq W_0 H^1(C_t) \subseteq W_1 H^1(C_t) \subseteq W_2 H^1(C_t) =
H^1(C_t).
$$
The defining properties of the monodromy weight filtration imply that
$$
W_0 H^1(C_t) = \im N \text{ and } W_1H^1(C_t) = \ker N.
$$
The Picard-Lefschetz formula implies that $W_0H^1(C_t)$ is spanned by the
Poincar\'e duals of the vanishing cycles.

Fix a tangent vector $\v$ of $0$ in $\D$ and denote the corresponding limit MHS
by $H^1(C_\v)$.

Denote the normalization of $C_0$ by $\Ctilde_0$. Let $D$ be the inverse image
in $\Ctilde_0$ of the double points of $C_0$. Results of the previous paragraph
imply that the natural morphisms
$$
\begin{CD}
H^1(C_0) @>\beta>> H^1(C_\v) @>\alpha>> H^1(C_0-D)
\end{CD}
$$
are morphisms of MHS.

The fact that $F^0 H^1(C_\v) = H^1(C_\v,\C)$ and $F^2 H^1(C_\v) = 0$ imply that
$$
\Gr^W_m H^1(C_\v)
$$
is of type $(0,0)$ when $m=0$, of type $(1,1)$ when $m=2$, and is a polarized
Hodge structure of with Hodge numbers $(1,0)$ and $(0,1)$ when $m=1$. Basic
topology implies that $\beta$ injective, from which it follows that $\beta$
induces a natural isomorphism
$$
\Gr^W_1 H^1(C_\v) \cong H^1(\Ctilde_0) \cong
\bigoplus_{\substack{\text{components}\cr E \text{ of }C_0}} H^1(E)
$$
of polarized Hodge structures. It also implies that $W_0 H^1(C_\v)$ is subspace
of $H^1(C_\v)$ generated by the Poincar\'e duals of the vanishing cycles and
that $\beta$ induces an isomorphism
\begin{equation}
\label{isom}
W_1 H^1(C_\v) \cong H^1(C_0)
\end{equation}

The classical Torelli Theorem, the fact that the theta divisor of the jacobian
of a smooth curve is irreducible and the semi-simplicity of polarized Hodge
structures of weight $1$ imply that the polarized limit MHS, $H^1(C_\v)$,
determines the normalization of $C_0$. Using this, and by approximating the
period mapping by the nilpotent orbit, one can give proofs of weak versions
(i.e., mod $t$) of the two results \cite[Cor.~3.2, Cor~3.8]{fay} in Fay's book.

If the normalization of $C_0$ is connected, then one can determine the divisor
classes in $\Jac \Ctilde_0$ of the pairs of points that are identified to
obtain $C_0$. This information is extracted from the extension
$$
0 \to W_0 H^1(C_\v) \to W_1H^1(C_\v) \to \Gr^W_1 H^1(C_\v) \to 0.
$$
using the work of Carlson \cite{carlson} and the isomorphism (\ref{isom}).
This implies that when $C_0$ is irreducible, the limit MHS $H^1(C_\v)$
determines $C_0$ up to isomorphism.

On the other hand, if the dual graph of $C_0$ is a tree, then all vanishing
cycles are trivial in homology, and the Picard-Lefschetz formula implies that
$N=0$. In this case, the monodromy weight filtration is trivial, and the limit
MHS $H^1(C_\v)$ is a polarized Hodge structure of weight 1. By the argument
above, this determines the non-rational components of $\Ctilde_0$  up to
isomorphism. But since there is no extension data, it gives no information
about how to reassemble $C_0$ from its irreducible components. In this case,
one can use the limit MHS on the truncation of the group ring of the
fundamental group of the smooth fiber by the $4$th power of its augmentation
ideal\footnote{The mixed Hodge structure on the fundamental groups is discussed
in the next section. Also, to apply it, one needs to choose a section $\sigma:
\D \to C$ of base points.} to determine $C_0$ up to finite ambiguity using the
limit mixed Hodge structure on the fundamental group which is constructed in
\cite{hain:dht2}. A proof of this assertion is not published, but most of the
technical constructions needed to prove it are worked out in detail in the
doctoral thesis of Rainer Kaenders \cite{kaenders}. More generally, I believe
that for all stable degenerations of curves, the limit mixed Hodge structure on
the truncation of the integral group ring of the fundamental group of the
smooth fiber by the 5th power of its augmentation ideal should determine $C_0$
up to finite ambiguity.

\section{Periods of $\pi_1(\Pminus,\b)$}

In this section, we sketch a fundamental example of computing the periods
of a limit MHS on the fundamental group of a the thrice punctured line.
The periods turn out to be Euler's the mixed zeta numbers (cf.\ \cite{zagier},
\cite{goncharov}).

\subsection{Hodge Theory of Homotopy Groups}

This is a very brief introduction to the Hodge theory of homotopy groups
of complex algebraic varieties.

First, suppose that $\pi$ is a discrete group. For a commutative ring $R$,
denote the group algebra of $\pi$ over $R$ by $R\pi$. There is a natural
augmentation
$$
\epsilon : R\pi \to R
$$
defined by taking $\sum r_j g_j$ to $\sum r_j$, where each $r_j \in R$ and
$g_j \in \pi$. The augmentation ideal $J_R$ is the kernel of the augmentation.
The mapping
$$
J_R/J_R^2 \to H_1(\pi,R)\quad (g-1) + J_R^2 \to [g]
$$
is a group isomorphism, which can be thought of as a kind of Hurewicz
isomorphism.

\begin{theorem}[Morgan \cite{morgan}, Hain \cite{hain:fund_gp}]
If $X$ is a complex algebraic variety and $x \in X$, then for each $s \ge 0$,
there is a natural MHS on the truncated group ring $\Z\pi_1(X,x)/J^{s+1}$.
These form an inverse system in the category of MHS. The mixed Hodge structure
on $J/J^2$ is dual to Deligne's MHS on $H^1(X)$.
\end{theorem}

These mixed Hodge structures are constructed using Chen's iterated integrals.

Suppose that $M$ is a manifold and that $w_1,\dots, w_r$ are smooth $\C$-valued
1-forms on $M$. For each piecewise smooth path $\gamma : [0,1] \to M$, we can
define
$$
\int_\gamma w_1\dots w_r =
\int \cdots \int_{0 \le t_1 \le \cdots \le t_r \le 1}
f_1(t_1) \cdots f_r(t_r)\, dt_1 \ldots dt_r,
$$
where $\gamma^\ast w_j = f_j(t)\, dt$ for each $j$. This is viewed as a
$\C$-valued function
$$
\int_\gamma w_1\dots w_r : PM \to \C
$$
on the the space of piecewise smooth paths in $M$. When $r=1$, $\int_\gamma w$
is just the usual line integral. An {\it iterated integral} is any function $PM
\to \C$ which is a linear combination of a constant function and basic iterated
integrals
$$
\int_\gamma w_1\dots w_r.
$$
An exposition of the construction of the MHS in $J/J^{s+1}$ can be found in
\cite{hain:bowdoin}.

\subsection{The MHS on $\pi_1(\Pminus,t)$}
\label{mhs-p1}

One can describe the MHS on the $J$-adic completion
$$
\Q\pi_1(\Pminus,t)\comp := \limproj{s} \Q\pi_1(\Pminus,t)/J^s.
$$
of $\pi_1(\Pminus,t)$ quite directly.

Consider the non-commutative power series ring
$$
A := \C\lb X_0, X_1 \rb
$$
freely generated by the indeterminates $X_0$ and $X_1$. Set
$$
w_0 = \frac{dz}{z}\text{ and } w_1 = \frac{dz}{1-z}.
$$
Now consider the $A$-valued iterated integral
$$
T = 1 + \int w_0 X_0 + \int w_1 X_1 + \cdots + \int
w_{j_1}\dots w_{j_r} X_{j_1} \dots X_{j_r} + \cdots
$$
It is not difficult to use the definition of iterated integrals to show that
the value $T(\gamma)$ of this on a path $\gamma$ depends only on its homotopy
class relative to its endpoints. In addition, one can show that if $\alpha$ and
$\beta$ are composable paths, then\footnote{Note that I use the topologists
convention that if $\alpha,\beta : [0,1] \to X$ are two paths with $\alpha(1) =
\beta(0)$, then $\alpha\beta$ is the path obtained by first traversing
$\alpha$, then $\beta$. This is the opposite of the convention used by many
papers in this field.}
\begin{equation}
\label{prod}
T(\alpha\beta) = T(\alpha)T(\beta).
\end{equation}

For each $t\in \Pminus$, we can thus define a homomorphism
$$
\pi_1(\Pminus,t) \to \text{ the group of units of $A$}
$$
by taking the class of the loop $\gamma$ based at $t$ to $T(\gamma)$. This
induces a homomorphism
$$
\Theta_t : \C\pi_1(\Pminus,t)\comp \to A
$$
which can easily be shown to be an isomorphism by using universal mapping
properties of free groups and free algebras.

We can use this to construct a MHS on $\Q\pi_1(\Pminus,t)\comp$. Give each
generator $X_j$ type $(-1,-1)$. Extend this bigrading to all monomials in $A$
in the standard way --- the monomial $X_I$ will have Hodge type $(-|I|,-|I|)$
where $I = (i_1,\dots,i_r)$ is a multi index and $|I|$ is its length $r$. Thus
$X_I$ has weight $-2|I|$. It is natural to define
$$
W_m A = \text{ the closure of the span of $X_I$ where $-2|I| \le m$}
$$
(which is the $\ell$th power of the maximal ideal when $m = -2\ell$) and
$$
F^p A = \text{ the span of the $X_I$ where $-|I| \ge p$}
$$
which is finite dimensional for all $p$. Pulling back these filtrations
along $\Theta_t$ defines the standard MHS on $\Q\pi_1(\Pminus,t)\comp$.

The periods of the MHS on $\Q\pi_1(\Pminus,t)\comp$ are the coefficients of the
monomials $X_I$ in the power series $T(\gamma)$, where $\gamma \in
\pi_1(\Pminus,t)$.

\begin{remark}
In the Hodge theory of complex varieties, the periods are typically (local)
coordinates in an appropriate moduli space of MHS, such as the one
$$
G_\Z \backslash G/F^0G
$$
discussed near the end of Section~\ref{monod_wt_filt}. However, in more
arithmetic situations, where the MHS arises from a variety defined (say) over
$\Q$, the de~Rham invariant together with its Hodge and weight filtrations has
a natural $\Q$-form. In this case, one usually takes the periods to be the
integrals over rational (betti) cycles of a basis of the $\Q$-de~Rham version
of the invariant which is adapted to the Hodge and weight filtrations.

This is the case above as $(\Pminus,t)$ is defined over $\Q$ when $t$ is
$\Q$-rational. There is an algebraic de~Rham theorem for the unipotent
fundamental group in this case (cf.\ \cite{hain:tianjin}), which is
particularly easy to describe in the case above. The dual of
$\C\pi_1(\P^1-\{0,1,\infty\},t)/J^{N+1}$ is, by Chen's $\pi_1$ de~Rham Theorem
(or elementary arguments in this case),
$$
V^\DR_\C :=
\bigg\{\sum_{r\le N} a_I \int w_{i_1}\dots w_{i_r} :
I \subset \{0,1\}^r,\ |I| = r,\ a_I \in \C\bigg\}.
$$
This has the natural $\Q$-form
$$
V^\DR_\Q :=
\bigg\{\sum_{r\le N} a_I \int w_{i_1}\dots w_{i_r} :
I \subset \{0,1\}^r,\ |I| = r,\ a_I \in \Q\bigg\}.
$$
The basis of $\C\pi_1(\P^1-\{0,1,\infty\},t)/J^{N+1}$ dual to basis $\{\int
w_{i_1}\dots w_{i_r}\}$ of $V^\DR_\Q$ is the set of monomials $X_I :=
X_{i_1}\dots X_{i_r}$. The periods of the MHS on
$\C\pi_1(\P^1-\{0,1,\infty\},t)/J^{N+1}$ are thus the integrals of the $\int
w_{i_1}\dots w_{i_r}$ over elements of
$\Q\pi_1(\P^1-\{0,1,\infty\},t)/J^{N+1}$.
\end{remark}

\subsection{The MHS on $\pi_1(\Pminus,\b)$}

In general, the periods of the MHS on $\pi_1(\Pminus,t)$ are difficult to
compute as the values of the coefficients of $T$ on paths based at $t$ do not
appear to be readily recognizable numbers when $t$ is general. However, if we
let $t$ approach 0, the periods (or their asymptotics) become more recognizable
as we shall explain.

When $t$ goes to 0, the MHS on $\Q\pi_1(\Pminus,t)\comp$ degenerates as all of
the integrals of length $\ge 2$ that begin or end with $w_0$, such as $\int w_0
w_1$, will diverge. (The MHS on the fundamental group can degenerate when the
variety becomes singular, and also when the base point runs off the edge of the
space as it is here.) Just as there is Schmid's theory of limits of Hodge
structures, there is a theory of limits of MHSs, at least in the geometric case
\cite{steenbrink-zucker}. (There is a fledgling theory in the abstract case
too, being developed by Kaplan and Pearlstein \cite{kaplan-pearlstein}.) We
will compute the periods of the limit mixed Hodge structure as $t\to 0$ with
respect to the tangent vector $\b := \partial/\partial t$, where $t$ is the
standard parameter on $\P^1$.

Deligne \cite{deligne:line} introduced the idea of the fundamental group of
(say) an affine curve $C$ with base point a non-zero tangent vector $\v$ at one
of the cusps $P$. It is simple and elegant. The fundamental group $\pi_1(C,\v)$
is the set of homotopy classes of loops $\gamma$ in $C \cup \{P\}$ based at
$P$, that leave $P$ with velocity vector $\v$ and return to $P$ with velocity
vector $-\v$. It is also required that $\gamma$ does not return to $P$ when
$0<t<1$. It is naturally isomorphic to the standard fundamental
group.\footnote{ An elegant way to prove this is to note that $\pi_1(C,\v)$ is
naturally isomorphic to $\pi_1(\Ctilde,[\v])$ where $\Ctilde$ is the real
oriented blowup of $C$ at $P$ and $[\v]$ is the point on the exceptional circle
$\partial \Ctilde$ corresponding to the oriented ray determined by $\v$.} This
will carry a limit MHS associated to the tangent vector $\v$.

\begin{definition}
For integers $n_1,\dots,n_r$, where $n_r>1$, define
$$
\zeta(n_1,\dots,n_r) =
\sum_{0<k_1<\dots <k_r} \frac{1}{k_1^{n_1}k_2^{n_2}\dots k_r^{n_r}}.
$$
\end{definition}

These numbers generalize the classical values of the Riemann zeta function at
integers $>1$ and were first considered by Euler. They have recently resurfaced
in the works of Zagier \cite{zagier} and Goncharov \cite{goncharov}. Their
$\Q$-linear span in $\R$ is a subalgebra $\MZN$ due to combinatorial identities
such as
$$
\zeta(n)\zeta(m) = \zeta(n,m) + \zeta(m,n) + \zeta(m+n).
$$
Since $\zeta(2) = \pi^2/6$,
$$
\MZN_\C := \MZN \oplus i\pi \MZN
$$
is a $\Q$-subalgebra of $\C$.

Mixed zeta numbers can be expressed as iterated integrals (cf.\
\cite{zagier}):
$$
\zeta(n_1,\dots,n_r) = \int_0^1
w_1 \overbrace{w_0 \dots w_0}^{n_1-1}
w_1 \overbrace{w_0 \dots w_0}^{n_2-1}
w_1 \dots
w_1 \overbrace{w_0 \dots w_0}^{n_r-1}.
$$
Here the path of integration is along the unit interval, and the integral
converges if and only if $n_r > 1$. This identity follows directly from the
definition of iterated integrals using the following slightly more general
form:

{\scriptsize
$$
\int_0^x
w_1 \overbrace{w_0 \dots w_0}^{n_1-1}
w_1 \overbrace{w_0 \dots w_0}^{n_2-1}
w_1 \dots
w_1 \overbrace{w_0 \dots w_0}^{n_r-1}
= \sum_{0<k_1<\dots <k_r} \frac{x^{k_r}}{k_1^{n_1}k_2^{n_2}\dots k_r^{n_r}}.
$$
}

The basic result we wish to explain is the following ``folk theorem.''
Goncharov in \cite{goncharov} has recently written down a proof of a
considerably more refined statement.

\begin{theorem}
\label{periods}
The space of periods of the limit mixed Hodge structure on
$\Q\pi_1(\Pminus,\b)\comp$ is precisely $\MZN_\C$.
\end{theorem}

It is convenient to give a proof of this theorem by appealing to the MHS on
spaces of paths. Denote the space of piecewise smooth paths from $a$ to $b$ in
$\Pminus$ by $P_{a,b} \Pminusp$. Endow it with the compact open topology. Then
$$
H_0(P_{a,b}\Pminusp,\Q)
$$
is the $\Q$-vector space generated by the homotopy classes of paths in
$\Pminus$ that go from $a$ to $b$. It is a module over $\Q\pi_1(\Pminus,a)$ and
can be completed in the $J$-adic topology. Denote this completion by
$$
H_0(P_{a,b}\Pminusp,\Q)\comp.
$$
(It is also a $\Q\pi_1(\Pminus,b)$ module. The completion with respect to this
action is identical with the completion above.)  As in the case where $a=b$, we
can define a MHS on it by pulling back the Hodge and weight filtrations of $A$
along the mapping
$$
\Theta_{a,b} : H_0(P_{a,b}\Pminusp,\Q)\comp \to A
$$
induced by taking $\gamma$ to $T(\gamma)$, which is an isomorphism. The natural
mappings
\begin{multline*}
H_0(P_{a,b}\Pminusp,\Q)\comp \otimes H_0(P_{b,c}\Pminusp,\Q)\comp \cr
\to H_0(P_{a,c}\Pminusp,\Q)\comp
\end{multline*}
are easily seen to be morphisms of MHS because of the formula (\ref{prod}).

By taking limit MHS, we can replace $a$ and $b$ by non-zero tangent vectors at
the cusps. In particular, we have a mixed Hodge structure on
$$
H_0(P_{\b,\c}\Pminusp,\Q)\comp
$$
where $\c$ is the tangent vector $\partial/\partial (1-t)$ at $t=1$.

The following result is essentially a restatement of a result of Le and
Murakami \cite[Thm.~A.9]{le-murakami}. (See also, \cite[Lemma~5.4]{goncharov}.)

\begin{lemma}
The periods of the path $[0,1]\in P_{\b,\c}\Pminusp$ under the mapping
$$
\Theta_{\b,\c} : H_0(P_{\b,\c} \Pminusp) \to A
$$
are mixed zeta numbers. All elements of $\MZN$ occur.
\end{lemma}

\begin{proof}
Let $\Phi \in A$ be the image of $[0,1]$ under the regularized period mapping
$\Theta_{\b,\c} : H_0(P_{\b,\c} \Pminus) \to A$. I claim that this is
Drinfeld's associator \cite{drinfeld}. The result will then follow as Le and
Murakami have shown that the coefficients  of the Drinfeld associator are
mixed zeta numbers and all mixed zeta numbers occur.

To see that $\Phi$ is the Drinfeld associator, we use the prescription given in
Corollary~\ref{limit} for computing the regularized limit periods of a flat
section. Applying this directly, we see that the renormalized value of
$T([0,1])$ is
$$
\lim_{t \to 0} t^{X_0} T([t,1-t]) t^{X_1}
$$
as the residue at $t=0$ of the connection for which $T[t,1-t]$ is a flat
section is
$$
\text{left multiplication by }X_0 + \text{ right multiplication by }X_1
$$
The limit is a well-known expression for Drinfeld's associator. (Cf.\
\cite{drinfeld}, \cite{le-murakami}.)
\end{proof}

\begin{center}
\begin{figure}[!ht]
\epsfig{file=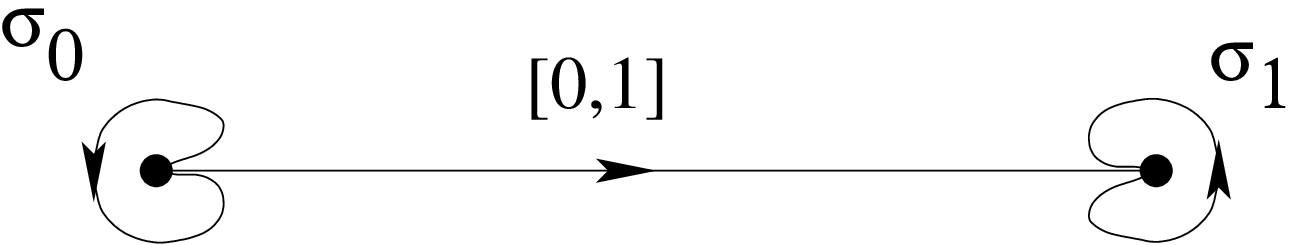, width=3in}
\caption{}\label{paths}
\end{figure}
\end{center}

Consider the paths illustrated in figure~\ref{paths}
where $\sigma_0 \in P_{\b,\b}\Pminus$ and $\sigma_1 \in P_{\c,\c}\Pminus$.

Theorem~\ref{periods} now follows as one can check, as in the proof of the
previous lemma, that under the renormalized homomorphism
$$
\Theta_\b : \pi_1(\Pminus,\b) \to A
$$
the image of the loop $\sigma_0$ is $\exp (2\pi i X_0)$ and that the image of
$\sigma_1$ under
$$
\Theta_\c : \pi_1(\Pminus,\c) \to A
$$
is $\exp(-2\pi i X_1)$. Since $\pi_1(\Pminus,\b)$ is generated by $\sigma_0$
and $[0,1]\sigma_1[0,1]^{-1}$, and since
$$
T([0,1]\sigma_1[0,1]^{-1}) = \Phi e^{2\pi i X_1} \Phi^{-1}
$$
we see that all periods of $\Q\pi_1(\Pminus,\b)\comp$ lie in $\MZN_\C$.

\section{Mixed Tate Motives}

The question arises as to why one should expect the periods of
$\pi_1(\Pminus,\b)$ to be so special and to have such a nice description in
terms of mixed zeta numbers. Deligne and Goncharov \cite{deligne-goncharov},
building on previous work of both authors, especially \cite{deligne:line} and
several of Goncharov's papers in {\sf arXiv.org}, have developed a satisfying
and deep explanation in terms of the theory of mixed Tate motives.

Thanks to the work of Voevodsky \cite{voevodsky}, Levine \cite{levine}, and
Deligne-Goncharov \cite{deligne-goncharov}, we now know that there is a
tannakian category of mixed Tate motives over $\Spec$ of the ring of
$S$-integers in a number field with the expected ext groups. So, in particular,
there is a tannakian category of mixed Tate motives over $\Spec \Z$. Such mixed
Tate motives have a Hodge realization, which is a mixed Hodge structure, all
of whose weight graded quotients are of type $(p,p)$. These mixed Hodge
structures will have two $\Q$-structures --- one coming from a $\Q$-de~Rham
structure, and another from topology (the ``Betti realization''). The entries
of the matrices relating them are, by definition, the periods of the motive.

Deligne and Goncharov \cite{deligne-goncharov} show that
$\Q\pi_1(\Pminus,\b)\comp$, equivalently the Lie algebra of the unipotent
fundamental group of $(\Pminus,\b)$, is a pro-object in the category of mixed
Tate motives over $\Spec \Z$. One reason one should suspect this is that
$$
\PminusZ := \Spec\Z[t,t^{-1},(1-t)^{-1}]
$$
has everywhere good reduction. However, the pointed variety $(\PminusZ,N)$,
where $N\in \Z -\{0,1\}$, will not have good reduction at primes dividing $N$
as then the base point will move to a cusp. In order for $(\Pminus,t)$ to have
everywhere good reduction, we are forced to take $t$ to be $\b$, or one of its
5 other images under $\Aut (\P^1,\{0,1,\infty\})$.

The fact that $\Q\pi_1(\Pminus,\b)\comp$ is a mixed Tate motive over $\Spec \Z$
has a remarkable consequence for the transcendence properties of mixed zeta
numbers \cite{terasoma} and for the Galois action on its $\Ql$-form
\cite{hain-matsumoto:gal}. In the first case, Zagier conjectured (unpublished),
and Terasoma \cite{terasoma} has proved, that the dimension of the $\Q$-linear
span of the mixed zeta numbers $\MZN_m$ of weight\footnote{The weight of
$\zeta(n_1,\dots,n_r)$ is, by definition, $n_1 + \dots + n_r$.} $m$ is bounded
by the dimension of the $m$th weight graded quotient of the graded algebra
$$
\Q[Z_2]\otimes \Q\langle Z_3,Z_5,Z_7,\dots \rangle
$$
where $Z_m$ has weight $m$ and $\Q\langle Z_3,Z_5,Z_7,\dots \rangle$ denotes
the free associative algebra generated by the $Z_\mathrm{odd}$.

One currently unresolved question is whether the periods of all mixed Tate
motives over $\Spec\Z$ lie in $\MZN_\C$. Goncharov
\cite[p.~385]{goncharov:periods} has conjectured that this is the case. To see
why this might be, we need to discuss the motivic Lie algebra of $\Spec\Z$ and
its connection with the absolute Galois group $G_\Q := \Gal \Qbar/\Q$ and its
action on $\Ql\pi_1(\Pminus,\b)\comp$.

Since the category of mixed Tate motives over $\Spec \Z$ is tannakian, it is
the category of representations of a pro-algebraic group $\A$ over $\Q$. Since
the tannakian fundamental group of the category of pure Tate motives over
$\Spec\Z$ is $\Gm$, this group is an extension 
$$
1 \to \U \to \A \to \Gm \to 1
$$
where $\U$ is a prounipotent group. That the ext groups in the category of
mixed Tate motives over $\Z$ have the desired relationship to $K_\dot(\Z)$
implies that the Lie algebra of $\U$ is a free pronilpotent Lie algebra on
generators $z_3,z_5,z_7,\dots$, where $z_n$ has weight $-2n$, cf.\
\cite{deligne-goncharov} and \cite{hain-matsumoto:survey}. The Lie algebra $\a$
of $\A$ is itself a pro-object of the category of mixed Tate motives, and
formal arguments show that all periods of mixed Tate motives over $\Spec \Z$
occur as periods of $\a$. Thus one should try to prove that all periods of the
MHS of $\a$ lie in $\MZN_\C$. One way to approach this is by studying the
action of the absolute Galois group on the algebraic fundamental group of the
thrice punctured line over $\Qbar$.

Every mixed Tate motive $V$ has an $\ell$-adic realization for each rational
prime $\ell$. This is a representation $\rho_\ell : G_\Q \to \Aut V_\Ql$, where
$V_\Ql$ denotes $V\otimes_\Q \Ql$. These actions induce a Zariski dense
homomorphism
$$
G_\Q \to \A(\Ql)
$$
through which $\A\otimes\Ql \to \Aut V$ factors. The homomorphism $\rho_\ell$
is the composite of this with the canonical homomorphism $\A \to \Aut V$:
$$
G_\Q \to \A(\Ql) \to \Aut V_\Ql.
$$
The image of $\A\otimes\Ql$ in $\Aut V_\Ql$ is thus the Zariski closure of the
image of $G_\Q$  (cf.\ \cite{hain-matsumoto:survey}).

In the case of the unipotent fundamental group of $(\Pminus,\b)$, $\rho_\ell$
is the natural Galois action
$$
G_\Q \to \A(\Ql) \to \Aut \Ql\pi_1(\Pminus,\b)\comp
$$
induced by the action of $G_\Q$ on its algebraic fundamental group. The Zariski
closure of the image of this action is the image of the homomorphism
$$
\A\otimes\Ql \to \Aut \Ql\pi_1(\Pminus,\b)\comp.
$$
This homomorphism is injective if and only if its derivative
$$
\a \to \Der\Ql\pi_1(\Pminus,\b)\comp
$$
is. Since the periods of $\Der\Ql\pi_1(\Pminus,\b)\comp$ lie in $\MZN_\C$, the
periods of $\a$ will too if the derivative is injective. At present it is not
known  whether or not this is the case, although there are computer results
that show that in ``small weights'' the derivative is injective.

\end{document}